\documentclass[amstex,12pt]{article}

\usepackage{graphicx}
\usepackage{amsfonts}
\usepackage{amssymb}
\usepackage{amsmath}
\usepackage{latexsym}
\usepackage{enumerate}

\newtheorem{thm}{Theorem}[section]
\newtheorem{lm}[thm]{Lemma}

\newtheorem{q}[thm]{Question}

\newtheorem{ex}[thm]{Example}
\newtheorem{df}[thm]{Definition}

\begin{document}

\author{Jin-ichi Itoh and Costin V\^\i lcu}
\title{On the number of cut locus structures \\on graphs}
\maketitle

\noindent {\bf Abstract.} {\small
We proved in \cite{iv2} that every connected graph can be realized as the cut locus 
of some point on some riemannian surface.
Here we give upper bounds on the number of such realizations.
\\{\bf Math. Subj. Classification (2000):} 53C22, 05C10}


\section{Introduction}

In this note, by a surface we always mean a complete, compact and connected $2$-dimensional riemannian manifold 
without boundary.
All the graphs we consider are undirected, connected, and may have multiple edges and loops, 
but not vertices of degree two unless explicitly stated otherwise.

The notion of cut locus was introduced by H. Poincar\'e \cite{p} in 1905, and gain since then an important place in global riemannian geometry.
The {\it cut locus $C(x)$} of the point $x$ in the riemannian manifold $M$ is the set of all extremities (different from $x$) of maximal (with respect to inclusion) shortest paths starting at $x$; for basic properties and equivalent definitions refer, for example, to \cite{ko} or \cite{sa}.

For riemannian surfaces $S$ is known that $C(x)$, if not a single point, is a local tree (i.e., each of its points $z$ has a neighborhood $V$ in $S$ such that the component $K_z(V)$ of $z$ in $C(x)\cap V$ is a tree), 
even a tree if $S$ is homeomorphic to the sphere.
(A {\it tree} is a set $T$ any two points of which can be joined
by a unique Jordan arc included in $T$.)

S. B. Myers \cite{M} for $d=2$, and M. Buchner \cite{Bu2} for general $d$, 
established that the cut locus of a real analytic riemannian manifold of dimension $d$ is homeomorphic to a finite $(d-1)$-dimensional simplicial complex.

We proved in \cite{iv2} a partial converse to Myers' theorem, namely that
{\sl every graph $G$ can be realized as a cut locus}; 
i.e., there exist a riemannian surface $S_G=(S_G,h)$ a point $x$ in $S_G$
such that $C(x)$ is isometric to $G$. The surface $S_G$ is in general not unique.
If, moreover, $G$ is cyclic of order $k$, then it can be realized as a cut locus on a surface of constant curvature.

Consider, for example, the graph represented in Figure 1 a).
It can be realized in two different ways as a cut locus, on an orientable, and on a 
non-orientable flat surface, see Figure 1 b) and c).

\begin{figure*}
\label{2_cycles}
\centering
  \includegraphics[width=1\textwidth]{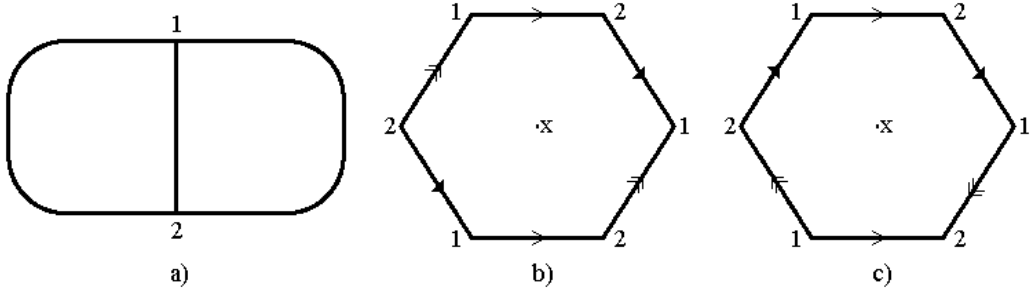}
\caption{A graph $G$ a), and two different realizations as a cut locus b), c).
The arrows show how to identify the edges of the regular hexagon such that the cut locus of its centre $x$ is $G$.}
\end{figure*}

In this note we obtain upper bounds on the number of distinct realizations of a graph as a cut locus, on general surfaces
(Section \ref{nb_CL}), and on orientable surfaces (Section \ref{nb_orient_CL}).
In particular, for the graph in Figure 1 a) there exist only two realizations, see Figure 1 b) and c),
up to graph isomorphisms.

The meaning of the words ``distinct'' or ``different'' above can be made precise by the use of 
{\sl patches and strips on graphs}, that we introduced in \cite{iv1}. 
In Section \ref{CL-str} we briefly present these notions and their properties, 
while in Section \ref{graph_prelim} we construct a set of paths on the extended intersection graph of cycles.
Both of these sections are preliminaries for Sections \ref{nb_CL} and \ref{nb_orient_CL}.

In a companion paper \cite{iv3} we are concerned about the orientability of the surfaces $S_G$ 
realizing the graph $G$ as a cut locus, and provide several criteria to recognize 
graphs having at least one such realization.


\section{Paths of cycles}
\label{graph_prelim}

In this section we recall a few --necessary-- definitions and facts about graphs, in order to fix notation.
The main goal is to construct, for any graph $G$, a set ${\cal P} \left({\cal G} \right)$ that
we shall use in the following sections.

\bigskip

{\bf Notation.} 
We denote by $V=V(G)$ the vertex set of the graph $G$, and by $E=E(G)$ the edge set of $G$.
We shall assume throughout the paper that
$m(G)$ denotes the number of edges in $G$, and $n(G)$ the number of vertices.

An edge in the graph $G$ is called {\it external} if it is incident to (least) one vertex of degree one, and is called a {\it bridge} if its removal disconnects $G$.
Denote by $b(G)$ the number of bridges in $G$.

A {\it path} in $G$ is a sequence $\left( e_{i_1},...,e_{i_k} \right)_{k \geq 1}$ of edges in $G$, 
any two consecutive edges of which share a common vertex. A {\it constant path} consists of a single vertex.


\bigskip

{\bf Edge contractions, cyclic graphs.}
An {\it edge contraction} in the graph $G$ is an operation which removes an edge from $G$ while simultaneously merging together the two vertices it used to connect to a new vertex; all other edges incident to either of the two vertices become incident to the new vertex.

The {\it cyclic part} of the graph $G$ is the minimal (with respect to inclusion) subgraph $G^{cp}$ of $G$, to which $G$ is contractible; i.e., the minimal subgraph of $G$ obtained by repeatedly contracting external edges, 
and for each vertex remaining of degree two (if any) contracting one of its incident edges.
To ease notation, we define
\begin{equation*}
m_{bc} (G) = m(G^{cp}) - b(G).
\end{equation*}

A graph is called {\it cyclic} if it is equal to its cyclic part.


\bigskip

{\bf Cycle space.}
The power set ${\cal E}$ of $E$ becomes a $Z_2$-vector space over the two-element field $Z_2$, 
if endowed with the symmetric difference as addition, and it is called the {\it edge space} of $G$. 

The {\it cycle space} is the subspace ${\cal Q}$ of ${\cal E}$ generated by (the edge sets of) all the simple cycles of $G$.
If $G$ is seen as a simplicial complex, ${\cal Q}$ is the space of $1$-cycles of $G$ with mod $2$ coefficients.
The symmetric difference $\ast$ of two simple cycles is either a simple cycle or a union of edge-disjoint simple cycles. 

The dimension $q=q(G)$ of the cycle space of the graph $G$ is given by $q(G)=m(G)-n(G)+1$.


\bigskip

{\bf Paths in the extended intersection graph of generating cycles.}
We shall make use of the {\it extended intersection graph} $H= H( {\cal G})$ of 
a set ${\cal G}=\{C_1, ...,C_q\}$ of generating cycles in $G$ (or $G^{cp}$).
Formally, $H$ has a vertex $i$ for each $C_i$, and 
an edge joining the vertices $i$ and $j$ for every distinct edge in $C_i \cap C_j$;
moreover, for every vertex shared by $C_i$ and $C_j$ which is not adjacent to an edge in $C_i \cap C_j$
\footnote
{Such a vertex has degree at least four, and can be seen as the contraction of a tree with degree three ramification points.},
we add one more edge between $i$ and $j$ in $H$ ($i \neq j$).
Thus, $H$ may have multiple edges and {\sl degree two vertices}.

A path $\Gamma$ in $H$ is {\it simple} if every vertex in $H$ appears at most once in $\Gamma$.

Denote by ${\cal A}$ the set of all {\sl simple paths} $\Gamma$ in $H$, 
and by ${\cal P} \left({\cal G} \right)$ the set of tuples of disjoint paths $\Gamma \in {\cal A}$.
The number $p({\cal G})$ of elements in the set ${\cal P} \left({\cal G} \right)$ will be useful in Section \ref{nb_CL}.
Formally,
\begin{equation*}
\label{P}
{\cal P} \left( {\cal G} \right) = \left\{ \left(\Gamma_1,...,\Gamma_k \right) | \;  k \leq m(G), \;
\Gamma_i \in {\cal G}, \; \Gamma_i \cap \Gamma_j =\emptyset, \; i \neq j, \; i,j \leq k \right\},
\end{equation*}
\begin{equation*}
\label{p(G)}
p({\cal G})= | {\cal P}  \left( {\cal G} \right) |.
\end{equation*}

\begin{ex}
\label{H_line}
Let ${\cal G}=\{C_1, C_2\}$ be the system of generating cycles for the graph $G$ in Figure 1 a), 
with $C_1$ the left cycle and $C_2$ the right cycle in the figure.
The (extended) intersection graph of  ${\cal G}$ consists of the vertices $1,2$ and the edge $12$ joining them.
The set ${\cal P} \left({\cal G} \right)=\{ (1), (2), (12), (1,2) \}$ consists of 
the constant paths $(1), (2)$, the path $(12)$, and the pair of disjoint paths $(1,2)$, hence
$p({\cal G})=4$.
\end{ex}

The following lemma gives an estimate on the order of $p({\cal G})$.

\begin{lm}
We have 
$$p({\cal G}) \geq 2^q -1,$$
with equality if and only if the extended intersection graph $H({\cal G})$ consists of vertices only. 
\end{lm}

\noindent{\sl Proof:}
Assume ${\cal G}=\{C_1,...,C_q\}$ is a system of generating cycles for the cycle space of $G$.
Then ${\cal P} \left({\cal G} \right)$ contains all constant paths, all pairs of constant paths, ..., 
all $(q-1)$-tuples of constant paths, and all $q$-tuples of constant paths.
And it contains a non-constant path if and only if the extended intersection graph $H({\cal G})$ contains at least one edge.
The conclusion follows from the binomial identity.
\hfill $\Box$


\section{Patches and strips}
\label{CL-str}

In this section we present notions and results developed in \cite{iv1}, and adapt them for our goal.

\medskip

A $G$-{\sl patch} on the graph $G$ is a topological surface $P_G$ with boundary, containing (a graph isometric to) $G$ and contractible to $G$.

A $G$-{\sl strip} (or a strip on $G$), is a $G$-patch whose boundary is topologically a circle.

An {\sl elementary strip} is an edge-strip (arc-strip) or a point-strip; 
i.e., a strip defined by the graph with precisely one edge (arc) of different extremities, 
respectively by the graph consisting of one single vertex.

An {\sl elementary decomposition} of a $G$-patch $P_G$ is a decomposition of $P_G$ into elementary strips such that:
\\- each edge-strip corresponds to precisely one edge of $G$;
\\- each point-strip corresponds to precisely one vertex of $G$.

\medskip

We shall represent patches in the plane as follows:
\\- consider an elementary decomposition of the patch;
\\- represent in the plane each vertex-strip;
\\- join the vertex strips by edge-strips in ``three dimensions'' (according to the patch), and
\\- flatten these edge-strips in the plane with overlappings.

\medskip

The following result was proven in \cite{iv1} for strips; the similar proof will not be repeated here.

\begin{lm}
\label{equiv_patches}
Consider an edge $e$ in the cyclic part of the graph $G$, which is not a bridge.
Assume we modify the elementary $e$-strip of a patch $P$ on $G$, by detaching one of its extremities,
rotating it with $2k \pi$, and attaching it back ($k \in I\!\!N)$.
Then the resulting patch is homeomorphic to the initial one.

If $e$ is  a bridge in $G$, then detaching one of its extremities, 
rotating the strip with $k \pi$, and attaching it back, produces a patch homeomorphic to the initial one ($k \in I\!\!N$).
\end{lm}

\begin{figure}
\label{Equivalent CL-structures}
\centering
  \includegraphics[width=0.4\textwidth]{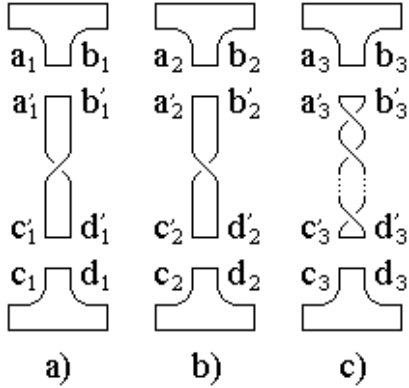}
\caption{Equivalent patches a), b) and c).
The edge-strip at a) corresponds to a rectangular band whose base is $\pi$-rotated ``to the left'' with respect to the top;
the edge-strip at b) corresponds to a rectangular band whose base is $\pi$-rotated ``to the right'' with respect to the top;
the edge-strip at c) corresponds to a rectangular band whose base is $(2k+1)\pi$-rotated ``to the left'' 
with respect to the top.}
\end{figure}

\begin{df}
The operation of detaching one extremity of an edge-strip from a patch,
rotating it with $\pi$ (or, equivalently --by Lemma \ref{equiv_patches}-- , with $(2k+1)\pi$, $k \in I\!\!N$), 
and attaching it back, will be called a {\it switch}, see Figure 2 (adapted from {\rm \cite{iv1}}).
\end{df}

{\sl We shall represent a switched edge-strip by drawing an ``{\rm x}'' over the edge, 
and a non-switched edge-strip by drawing an ``{\rm =}'' over the edge, see Figure 3.}

\bigskip

We explain now the relationship between patches and cut locus realizations of graphs.
Assume first that the cut locus $C(x)$ of the point $x$ in the surface $S$ is isometric to the graph $G$.
Then, cutting off the surface an open intrinsic disc of radius smaller than the injectivity radius at $x$,
provides a strip on $G$. The converse is established by the following result.

\begin{thm}
\label{glue} {\rm \cite{iv2}}
For every graph $G$ there exists at least one $G$-strip, and 
each $G$-strip provides a realization of $G$ as a cut locus.
\end{thm}

\begin{ex}
With the convention above, the strips defined by the surfaces illustrated in Figure 1 b) and c) 
are represented in Figure 3.
\end{ex}

\begin{figure*}
\label{2_cycles_CLS}
\centering
  \includegraphics[width=0.8\textwidth]{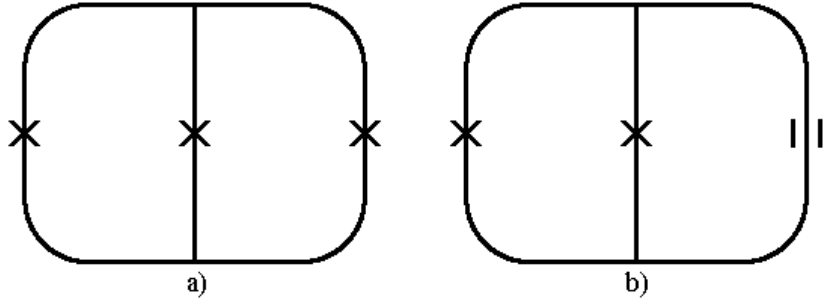}
\caption{Schematic representation for the strips defined by Figure 1 b) and c).}
\end{figure*}

Considering the planar representation of patches described above and Lemma \ref{equiv_patches}, 
counting patches on a given graph reduces to counting the odd or even number of switches of its edge-strips,
for non-bridge edges. Thus, we obtain the following.

\begin{lm}
\label{number_patches}
The number of distinct patches on the graph $G$ is $2^{m_{bc}(G)}$.
\end{lm}


\section{On the number of CL-structures}
\label{nb_CL}

In the view of Theorem \ref{glue} and Lemma \ref{equiv_patches}, 
counting different realizations of the graph $G$ as a cut locus reduces to counting different strips on $G^{cp}$. 

We call each $G^{cp}$-strip a {\it cut locus structure} (shortly, a {\it CL-structure}) on $G$; 
see \cite{iv1} for basic properties of such structures.

We need one more preliminary result.

\begin{lm}
\label{cycle_path}
Let ${\cal G}=\{C_1, ...,C_q\}$ be a set of generating cycles in the graph $G$, and
${\cal P} \left( {\cal G} \right)$ the set defined at the end of Section \ref{graph_prelim}.
For any $(\Gamma_1, ...,\Gamma_k)$ in ${\cal P} \left( {\cal G} \right)$, 
there exists a patch $P$ on $G$ with the following properties:
\\(i) $P$ has precisely $k+1$ boundary components, and
\\(ii) if the path $\Gamma_j$ passes through the vertices $p_1, ..., p_j$ (in this order) 
of the extended intersection graph of ${\cal G}$, 
then there is a boundary component of $P$ which encloses the cycle 
$C=C_{p_1} \ast C_{p_2} \ast ...\ast C_{p_j}$ (with $k \leq q$, $j \leq k$, and $\{p_1,...,p_j\} \subset \{1,...,q\}$).
\end{lm}

\noindent{\sl Proof:}
Contract all edges of $G$ in $\Gamma=\Gamma_1 \cup ...\cup \Gamma_k$, and denote by $\bar G$ the resulting graph.
By Theorem \ref{glue}, there exists a strip $S_{\bar G}$ on $\bar G$.
Define now a patch $P$ on $G$ by the use of $S_{\bar G}$: it coincides with $S_{\bar G}$ on the non-contracted edges of $G$,
and it is completed by non-switched edge-strips along the edges in $\Gamma$. 
It is easily seen that this completion produces a boundary component for every path $\Gamma_j$.
\hfill $\Box$

\begin{thm}
\label{number_CLS}
An upper bound for the number of distinct cut locus structures on a graph $G$ is given by
$$2^{m_{bc}(G)} - p({\cal G}),$$
and this is sharp.
\end{thm}

\noindent{\sl Proof:}
For the first part, notice that each patch produced by an element 
$(\Gamma_1, ...,\Gamma_k) \in {\cal P} \left( {\cal G} \right)$ as in Lemma \ref{cycle_path} is not a strip,
because it has at least two boundary components.
Any two such patches are different, because the tuples of paths producing them are different.
We obtain the desired inequality by the use of Lemma \ref{number_patches}.

For the last part of the statement, we consider again the graph $G_0$ in Figure 1 a).
From Example \ref{H_line}, we obtain for it the upper bound $4=2^3 - 4$.
On the other hand, the CL-structure in Figure 1 b) corresponds to switched edge-strips on $G_0$,
while the CL-structure in Figure 1 c) corresponds to two edge-strips on $G_0$ switched, and one edge-strip non-switched,
see Figure 3.
Since this last edge-strip can be chosen in three different ways, there exist (one+three=) four strips on $G_0$.
(They actually reduce to only two, if considered up to graph isomorphisms).
\hfill $\Box$


\section{On the number of orientable CL-structures}
\label{nb_orient_CL}

In this section we obtain an upper bound for the number of {\it orientable cut locus structures} on the graph $G$: 
$G^{cp}$-strips which are themselves orientable surfaces. Consequently \cite{iv2}, they produce realizations of $G$
as a cut locus on orientable riemannian surfaces.

\begin{lm}
\label{switch_edge}
Let ${\cal C}$ be a CL-structure on the graph $G$, and $E$ a non-switched edge-strip in ${\cal C}$. 
Switching $E$ yields a new CL-structure ${\cal C}'$ on $G$, 
which is non-orientable if ${\cal C}$ is orientable.
\end{lm}

\begin{figure*}
\label{Orient_vs_Non-orient}
\centering
  \includegraphics[width=0.8\textwidth]{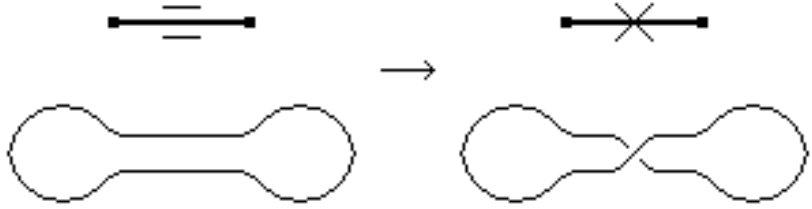}
\caption{switching a non-switched edge-strip in a CL-structure.}
\end{figure*}

\noindent{\sl Proof:}
In the left part of Figure 4 is represented 
the edge-strip $E$ as part of the boundary of ${\cal C}$.
It is easy to see that the resulting patch after the switch of $E$ is still a strip, see the right part of Figure 4.

The non-orientability of ${\cal C}'$ follows now in two ways.

Roughly, one can simply notice that the resulting surface has one face, if the initial surface has two faces.

Formally, one considers a realization of $G$ as a cut locus obtained via ${\cal C}$; 
i.e., a riemannian surface $S_G$ and a point $x$ in $S_G$ such that $C(x)=G$.
By our assumption, $S_G$ is orientable, hence the edge $e$ of $G$ (corresponding to $E$) 
has two images in the tangential cut locus (in the unit tangent space at $x$), 
with different orientations.
After the switch, the orientation of one image changes, yielding a non-orientable surface.
(Figure 1 b) and c) illustrates the effect of a switch for the graph in Figure 1 a).)
\hfill $\Box$

\bigskip

Lemma \ref{switch_edge} suggests the following definition and open question.

\begin{df}
Let ${\cal C}$ be a CL-structure on the graph $G$, $E$ a switched edge in ${\cal C}$,
and ${\cal C}'$ the patch on $G$ obtained from ${\cal C}$ by un-switching $E$. 
Then ${\cal C}'$ has at most two boundary components.

The edge $e$ is called {\sl transversal} to ${\cal C}$ if ${\cal C}'$ is not a strip on $G$, 
and it is called {\sl longitudinal} to ${\cal C}$ if ${\cal C}'$ is a strip on $G$.
\end{df}

\begin{q}
Do all orientable CL-structures have longitudinal edges?
If yes, estimate their number. 
\end{q}

\begin{ex}
The non-orientable CL-structure in Figure 3 b) has two switched edge-strips, both of them transversal.
The orientable CL-structure in Figure 3 a) has three switched edge-strips, all of them longitudinal.
\end{ex}

Denote by $P_{\rm x}=P_{\rm x}(G)$ the patch on the graph $G$ having all edge-strips switched,
by ${\cal O}={\cal O} (G)$ the set of all orientable $G$-strips, any
by ${\cal N}={\cal N} (G)$ the set of all non-orientable $G$-strips.

\begin{lm}
\label{Phi}
Consider a graph $G$, and the mapping
$\Phi : {\cal O} \setminus \{P_{\rm x}\} \to {\cal N}$ defined by Lemma \ref{switch_edge}. 
Assume $\Phi ({\cal C}_1)=\Phi ({\cal C}_2)$, with ${\cal C}_1 \neq {\cal C}_2 \in {\cal O} \setminus \{P_{\rm x}\}$, 
and let $e_1$ be the edge switched by $\Phi$ in ${\cal C}_1$, and $e_2$ the edge switched by $\Phi$ in ${\cal C}_2$. Then:

a)  ${\cal C}_1$ and ${\cal C}_2$ coincide outside $e_1 \cup e_2$;

b) $e_1$ and $e_2$ can be included in a simple cycle of $G$;

c) removing $e_1$ and $e_2$ from $G$ disconnects $G$.
\end{lm}

\noindent{\sl Proof:}
a) Take ${\cal C} \in {\rm Im}(\Phi) \subset {\cal N}$, such that $\Phi ({\cal C}_1)=\Phi ({\cal C}_2)={\cal C}$.
Since ${\cal C}$ and ${\cal C}_i$ differ only along the edge $e_i$, $i=1,2$, 
${\cal C}_1$ and ${\cal C}_2$ coincide outside $e_1 \cup e_2$.

b) Let $C$ be a cycle of $G$ containing the edge $e_1$;
it exists, because $e_1$ is not a bridge of $G$.
If $C$ doesn't contain $e_2$ then the cycle patch $P^1_C$ in ${\cal C}_1$ is an orientable surface,
while the cycle patch $P_C$ in ${\cal C}$ is non-orientable. 
But ${\cal C}_2$ coincides with ${\cal C}$ on $P_C$ (because they differ only along $e_2$),
hence $P_C$ is included in ${\cal C}_2$ which is orientable, and a contradiction is obtained.

c) Let $v^i$ and $w^i$ be the vertices of $G$ joined by $e_i$, $i=1,2$; then $v^i \neq w^i$, because of b).
If removing $e_1$ and $e_2$ from $G$ doesn't disconnect $G$, 
there exists a path $\Gamma$ in $G$ disjoint to $e_1 \cup e_2$ and joining $v^1$ to $w^1$, or $v^2$ to $w^2$.
Assume the first case holds; then $\Gamma \cup e_1$ is a cycle in $G$ disjoint to $e_2$.
A contradiction is now obtained, just as the one proving b).
\hfill $\Box$

\bigskip

Let $c({\cal G})$ be the maximal length (=number of edges) of a cycle in the system of simple generating cycles ${\cal G}$.
By $\lfloor \alpha \rfloor$ we denote the largest integer smaller than $\alpha \in I \!\! R$.

\begin{lm}
\label{Phi-1}
Consider a graph $G$, a system of simple generating cycles ${\cal G}$, and the mapping
$\Phi : {\cal O} \setminus \{P_{\rm x}\} \to {\cal N}$ defined by Lemma \ref{switch_edge}. 
Then $\Phi^{-1} ({\cal C})$ contains at most $c({\cal G})$ elements, for any ${\cal C} \in {\cal N}$.

If, moreover, $G$ is a $3$-graph without bridges then 
$\Phi^{-1} ({\cal C})$ contains at most $\lfloor \frac{c({\cal G})}2 \rfloor$ elements, for any ${\cal C} \in {\cal N}$.
\end{lm}

\noindent{\sl Proof:}
Assume there are distinct CL-structures ${\cal C}_i \in {\cal O} \setminus \{P_{\rm x}\}$
such that $\Phi ({\cal C}_i)={\cal C}$, and let $e_i$ be the edge switched by $\Phi$ in ${\cal C}_i$ ($i=1,..., k$).
Then, reasoning as in the proof of Lemma \ref{Phi}, we obtain that the edges $e_i$ are either all, or none of them, 
included in a cycle $C$ of $G$. Therefore, $k \leq c({\cal G})$.

Assume now that $G$ is a $3$-graph without bridges.
We claim that the cycle $C$ contains at least one edge between any two edges $e_i$ as above, 
which would prove that $k \leq \lfloor \frac{c({\cal G})}2 \rfloor$.
Assume the contrary be true, hence there exist two adjacent edges as above, say $e_1$ and $e_2$.
Denote by $v_{12}$ their common extremity, and by $f_{12}$ the third edge of $G$ at $v_{12}$.
Since $f_{12}$ is not a bridge of $G$, there exits a path $\Gamma$ in $G$ not containing $f_{12}$, 
joining $v_{12}$ to the other vertex of $f_{12}$. Then $\Gamma$ contains one of $e_1, e_2$.
It follows that removing $e_1$ and $e_2$ would leave $G$ connected, a contradiction to Lemma \ref{Phi} c)
proving the claim.
\hfill $\Box$

\begin{ex}
The mapping $\Phi$ defined by Lemma \ref{switch_edge} is not injective.
Indeed, we prove in {\rm \cite{iv3}} that the CL-structures in Figure \ref{2_cycles_CLS} a) and b) are both orientable.
Their image under the mapping $\Phi$ is the non-orientable CL-structure in Figure \ref{2_cycles_CLS} c).
\end{ex}

\begin{figure*}
\label{Non-injective_Phi}
\centering
  \includegraphics[width=1.0\textwidth]{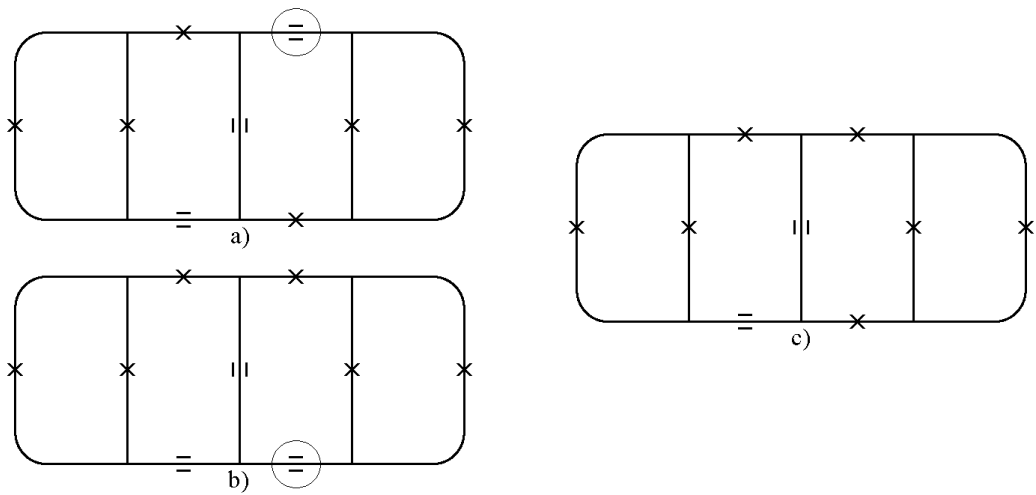}
\caption{Two orientable CL-structures, a) and b), both having the same image c) under the mapping $\Phi$.
The circles indicate in each case the edge to be switched.}
\end{figure*}

Denote by $O(G)$ the number of orientable CL-structures on the graph $G$, 
and by $N(G)$ the number of non-orientable CL-structures on $G$.

\bigskip

By Lemma \ref{Phi-1}, at most $c(G)$ orientable CL-structures on $G$ different from $P_x$ correspond to each non-orientable CL-structures on $G$. Thus, we immediately obtain the following.

\begin{lm}
\label{O_vs_N}
For every graph $G$ holds $O(G) \leq c({\cal G}) N(G) +1$.

If, moreover, $G$ is a $3$-graph without bridges then $O(G) \leq \lfloor \frac{c({\cal G})}2 \rfloor N(G) +1$.
\end{lm}

\begin{q}
In all examples we considered, the number of orientable CL-structures on a graph $G$ is 
(even ``much'') less than the number of non-orientable CL-structures on $G$. Is this always true?
\end{q}

It is clear that every CL-structure on a graph with an odd number $q$ of generating cycles is non-orientable.
We can now obtain an estimate for $O(G)$ if $q$ is even, by simply putting together 
Theorem \ref{number_CLS} and Lemma \ref{O_vs_N}.

\begin{thm}
\label{number_or-CLS}
An upper bound for the number of distinct orientable cut locus structures on 
a graph $G$ with an even number of generating cycles is
$$O(G) \leq \frac{c({\cal G}) \left( 2^{m_{bc}(G)} - p({\cal G}) \right)}{c({\cal G})+1} +1,$$
where ${\cal G}$ is a system of simple generating cycles in $G$.

If, moreover, $G$ is a $3$-graph without bridges then 
$$O(G) \leq \frac{\lfloor \frac{c({\cal G})}2 \rfloor \left( 2^{m_{bc}(G)} - p({\cal G}) \right)}
{\lfloor \frac{c({\cal G})}2 \rfloor+1} +1.$$
\end{thm}

We end with two open questions.

\begin{q}
How many orientable CL-structures can coexist on a Riemannian surface of genus $g$?
\end{q}

\begin{q}
How many (orientable) CL-structures exist on Riemannian surfaces of genus $g$?
\end{q}

In the companion paper \cite{iv3} we provide several criteria to recognize 
graphs having at least one orientable realization as a cut locus, and 
use them to provide the list of all orientable CL-structures on graphs
with $4$ generating cycles (which live as cut loci on surfaces of genus $2$).


\bigskip

\noindent {\bf Acknowledgement. } J. Itoh was partially supported by 
Grant-in-Aid for Scientific Research from the Japan Society for the Promotion of Science.

C. V\^\i lcu was partially supported by the 
grant PN II Idei 1187 of the Romanian Government.




\bigskip

Jin-ichi Itoh

\noindent {\small Faculty of Education, Kumamoto
University
\\Kumamoto 860-8555, JAPAN
\\j-itoh@gpo.kumamoto-u.ac.jp}

\medskip

Costin V\^\i lcu

\noindent {\small Institute of Mathematics ``Simion Stoilow'' of the
Romanian Academy
\\P.O. Box 1-764, Bucharest 014700, ROMANIA
\\Costin.Vilcu@imar.ro}

\end{document}